\documentclass[11pt,a4paper]{article}

\usepackage{inputenc}
\usepackage{amsmath}
\usepackage{amsfonts}

\usepackage{hyperref}

\title{Using Max-Algebra Linear Models\\ in the Representation of Queueing Systems\thanks{Proc. 5th SIAM Conf. on Applied Linear Algebra, Snowbird, UT, June 15-18, 1994, 155--160}
}

\author{Nikolai K. Krivulin\thanks{Postdoctoral Research Fellow, Faculty of Mathematics and Mechanics, St.Petersburg State University, Bibliotechnaya sq.2, Petrodvorets, St.Petersburg, 198904 Russia}}

\date{}

\def\sumo_#1^#2{\setbox0=\hbox{$\displaystyle{\sum}$}
                \setbox1=\hbox{$\scriptstyle{#1}$}
                \setbox2=\hbox{$\scriptstyle{#2}$}
		\setbox3=\hbox{${}_{{}_\oplus}\mathsurround=0pt$}
		\dimen1=.5\wd1 \advance\dimen1 by-.5\wd0
		\ifdim\dimen1>0pt
		   \ifdim\dimen1>\wd3 \kern\wd3 \else\kern\dimen1\fi\fi
		\dimen2=.5\wd2 \advance\dimen2 by-.5\wd0
		\ifdim\dimen2>0pt
		   \ifdim\dimen2>\wd3 \kern\wd3 \else\kern\dimen2\fi\fi
		\mathop{{\sum}{}_{{}_\oplus}}_{\kern-\wd3 #1}^{\kern-\wd3 #2}}
\def\sumol_#1{\setbox0=\hbox{$\displaystyle{\sum}$}
             \setbox1=\hbox{$\scriptstyle{#1}$}
	     \setbox3=\hbox{${}_{{}_\oplus}\mathsurround=0pt$}
	     \dimen1=.5\wd1 \advance\dimen1 by-.5\wd0
	     \ifdim\dimen1>0pt
	        \ifdim\dimen1>\wd3 \kern\wd3 \else\kern\dimen1\fi\fi
	     \mathop{{\sum}{}_{{}_\oplus}}_{\kern-\wd3 #1}}
\def\prodo_#1^#2{\setbox0=\hbox{$\displaystyle{\prod}$}
                \setbox1=\hbox{$\scriptstyle{#1}$}
                \setbox2=\hbox{$\scriptstyle{#2}$}
		\setbox3=\hbox{${}_{{}_\otimes}\mathsurround=0pt$}
		\dimen1=.5\wd1 \advance\dimen1 by-.5\wd0
		\ifdim\dimen1>0pt
		   \ifdim\dimen1>\wd3 \kern\wd3 \else\kern\dimen1\fi\fi
		\dimen2=.5\wd2 \advance\dimen2 by-.5\wd0
		\ifdim\dimen2>0pt
		   \ifdim\dimen2>\wd3 \kern\wd3 \else\kern\dimen2\fi\fi
	      \mathop{{\prod}{}_{{}_\otimes}}_{\kern-\wd3 #1}^{\kern-\wd3 #2}}

\begin{document}

\maketitle

\begin{abstract}
The application of the max-algebra to describe queueing systems by both linear
scalar and vector equations is discussed. It is shown that these equations may
be handled using ordinary algebraic manipulations. Examples of solving the
equations representing the $ G/G/1 $ queue and queues in tandem are also
presented.
\end{abstract}

\section{Introduction}
Max-algebra \cite{Cunn91,Olsd92} is the system
$ (\mathbb{R} \cup \{\varepsilon\}, \oplus, \otimes) $, where
$$
\varepsilon = -\infty, \;\;\; x \oplus y = \max(x,y),
\;\;\; x \otimes y = x + y \;\;\; \forall x,y \in \mathbb{R}.
$$
It has the following properties which can be easily verified

$$
\begin{array}{ll}
\forall x,y,z \in \mathbb{R} & x \oplus (y \oplus z) = (x \oplus y) \oplus z,
                                                                        \;\;\;
                                                 x \oplus y = y \oplus x,  \\
             & x \otimes (y \otimes z) = (x \otimes y) \otimes z, \;\;\;
                                                 x \otimes y = y \otimes x, \\
             & x \otimes (y \oplus z) = (x \otimes y) \oplus (x \otimes z), \\
	     & x \oplus \varepsilon =  x, \;\;\;  x \oplus x = x, \\
	     & x \otimes e = x, \;\; \mbox{where $ e=0 $}.
\end{array}
$$
In the max-algebra these properties allow ordinary algebraic manipulation of
linear expressions to be performed under the usual conventions regarding
brackets and precedence of $ \otimes $ over $ \oplus $. Moreover, the
scalar max-algebra is extended to the max-algebra of vectors in the regular
way. To emphasize parallels between conventional linear algebra and the
max-algebra, similar notations are used for the iterated operations
$ \oplus $ and $ \otimes $
$$
\sumo_{i=1}^{n} x_{i} = x_{1} \oplus \cdots \oplus x_{n},
\qquad
\prodo_{i=1}^{n} x_{i} = x_{1} \otimes \cdots \otimes x_{n}.
$$

The max-algebra theory is currently under investigation. There are a number of
classical algebraic results reformulated and proved in this algebra.
Specifically, the eigenvalue problem has been solved, and analogues of
Cramer's rule and the Cayley-Hamilton theorem have been found (see survey
papers \cite{Cunn91,Olsd92}). Moreover, as a research tool in studying
practical problems, the max-algebra finds expanding applications in many
fields of operations research and optimization, including the analysis and
performance evaluation of discrete event dynamic systems
\cite{Cunn91,Gree91,Olsd92}.

Although max-algebra models were successfully applied to investigate certain
classes of discrete event dynamic systems, the models of queues have received
little or no attention. The purpose of this paper is to show how queueing
systems may be described using the max-algebra approach by linear algebraic
equations. As illustrations of handling these algebraic models, the solutions
of the equations representing the $ G/G/1 $ queue and queues in tandem are
presented.

\section{A Linear Algebraic Model for the $ G/G/1 $ Queue}
We start with the linear max-algebra representation of the $ G/G/1 $
queue which provides the basis for more complicated models of queueing
systems. In the analysis of queueing systems, it is common to apply recursive
equations to describe their dynamics analytically. Such equations are normally
written in terms of recursions for the arrival and departure times of
customers, and involve the operations of maximum and addition
\cite{Gree91,Kriv90,Kriv93}.

To set up the equations that represent the $ G/G/1 $ queue in the ordinary
way, consider a single server queue with infinite buffer capacity. Once a
customer arrives into the system, he occupies the server provided that it is
free. If the customer finds the server busy, he is placed into the
buffer and has to wait until the service of all his predecessors completes.

Denote the interarrival time between the $k$th customer and his predecessor by
$ \alpha_{k} $, and the service time of the $k$th customer by
$ \tau_{k} $. Furthermore, let $ A(k) $ and $ D(k) $ be the arrival
and departure times of the $k$th customer, respectively. As is customary, we
assume that $ \alpha_{k} \geq 0 $ and $ \tau_{k} \geq 0 $ are given
parameters, whereas $ A(k) $ and $ D(k) $ are unknown variables. With
the conditions that the queue starts operating at time zero and it is free at
the initial time, one can readily represent the system dynamics by the set of
equations \cite{Gree91,Kriv90,Kriv93}
\begin{eqnarray*}
  A(k) & = & A(k-1) + \alpha_{k}, \\
  D(k) & = & \max(A(k), D(k-1)) + \tau_{k}.
\end{eqnarray*}

Let us now replace the usual operation symbols by those of the max-algebra and
rewrite the equations in their equivalent form as
\begin{eqnarray}
    A(k) & = & \alpha_{k} \otimes A(k-1), \label{set1-1} \\
    D(k) & = & \tau_{k} \otimes (A(k) \oplus D(k-1)). \label{set1-2}
\end{eqnarray}
Under the properties of the operation $ \oplus $ and $ \otimes $ these
equations could be handled much as if they were ordinary linear equations in
the conventional algebra. Specifically, applying a usual technique to solve
the equations for the unknown variables $ A(k) $ and $ D(k) $, we get
\begin{equation}\label{sol1}
A(k) = \alpha_{1} \otimes \cdots \otimes \alpha_{k}, \;\;\;
D(k) = \sumo_{i=1}^{k}  \alpha_{1} \otimes \cdots \otimes \alpha_{i}
                        \otimes \tau_{i} \otimes \cdots \otimes \tau_{k}.
\end{equation}

\subsection{The Matrix Representation}
To produce a matrix representation for the $ G/G/1 $ queue let us first
define the vector $ \mbox{\boldmath $D$}(k) = (D_{0}(k), D_{1}(k))^{T}$
with components $ D_{0}(k)=A(k)$, $ D_{1}(k)=D(k) $, and replace the symbols
$ \alpha_{k} $ and $ \tau_{k} $ by $ \tau_{0k} $ and
$ \tau_{1k} $, respectively, $ k=1,2, \ldots $. It is convenient to
preassign $ D_{0}(0) = D_{1}(0) = e $, and
$ D_{0}(k) = D_{1}(k) =\varepsilon $ for all $ k < 0 $. With the new
notations, the equations (\ref{set1-1}-\ref{set1-2}) may be rewritten as
\begin{eqnarray}
D_{0}(k) & = & \tau_{0k} \otimes D_{0}(k-1), \label{set11-1} \\
D_{1}(k) & = & \tau_{1k} \otimes (D_{0}(k) \oplus D_{1}(k-1)). \label{set11-2}
\end{eqnarray}
Substitution of (\ref{set11-1}) into (\ref{set11-2}) and the implementation of
distributivity of $ \otimes $ over $ \oplus $ give
\begin{eqnarray*}
D_{0}(k) & = & \tau_{0k} \otimes D_{0}(k-1), \\
D_{1}(k) & = & \tau_{1k} \otimes \tau_{0k} \otimes D_{0}(k-1)
                          \oplus \tau_{1k} \otimes D_{1}(k-1).
\end{eqnarray*}
We may now represent the model in matrix notations by the equation
\begin{equation}\label{equ1}
\mbox{\boldmath $D$}(k) = T_{k} \otimes \mbox{\boldmath $D$}(k-1),
\end{equation}
where the transition matrix is defined as
$$
T_{k} = \left( \begin{array}{rc}
		 \tau_{0k} & \varepsilon \\
		 \tau_{1k} \otimes \tau_{0k} & \tau_{1k}
	       \end{array} \right).
$$

\section{Linear Models of $ G/G/1 $ Queues in Tandem}
In this section we extend the max-algebra linear models to cover systems of
$ G/G/1 $ queues operating in tandem. As the basic system of this type,
we first consider a series of $ n $ queues with infinite buffers. Each
customer that arrives into this system is initially placed in the buffer at
the $1$st server and then has to pass through all the queues consecutively.
Upon the completion of his service at server $ i $, the customer is
instantaneously transferred to queue $ i+1 $, $ i=1, \ldots, n-1$. The
customer leaves the system after his service completion at the $n$th server.

For the tandem queueing system the equations (\ref{set11-1}-\ref{set11-2})
can be easily generalized as
\begin{eqnarray}
   D_{0}(k) & = & \tau_{0k} \otimes D_{0}(k-1), \label{set2-1} \\
   D_{i}(k) & = & \tau_{ik} \otimes (D_{i-1}(k) \oplus D_{i}(k-1)),
                                        \;\;   i=1, \ldots, n. \label{set2-2}
\end{eqnarray}
where $ D_{i}(k) $ and $ \tau_{ik} $ denote the departure time and the
service time of $k$th customer at server $ i $, respectively.

Let $ \mbox{\boldmath $D$}(k) = (D_{0}(k), \ldots, D_{n}(k))^{T}$ be the
vector of the $k$th departure times in the system. Similarly as in the case of
the $ G/G/1 $ queue, we may write the vector equation representing the
tandem queueing system in the form (\ref{equ1}) with the lower triangular
transition matrix
$$
T_{k} = \left( \begin{array}{rrrcc}
	       \tau_{0k} & \varepsilon & \varepsilon & \ldots & \varepsilon \\
	       \tau_{1k} \otimes \tau_{0k} & \tau_{1k} &
                                         \varepsilon & \ldots & \varepsilon \\
	       \vdots & \vdots &     & \ddots  & \vdots \\
	       \tau_{n-1k} \otimes \cdots \otimes \tau_{0k} &
	       \tau_{n-1k} \otimes \cdots \otimes \tau_{1k} &
               \tau_{n-1k} \otimes \cdots \otimes \tau_{2k} &
                                                            & \varepsilon \\
	       \tau_{nk} \otimes \cdots \otimes \tau_{0k} &
	       \tau_{nk} \otimes \cdots \otimes \tau_{1k} &
               \tau_{nk} \otimes \cdots \otimes \tau_{2k} &
                                                         \ldots & \tau_{nk}
	       \end{array} \right).
$$

Furthermore, we may find the solution of the set of recursive equations
(\ref{set2-1}-\ref{set2-2}) as an extension of (\ref{sol1}). With usual
algebraic manipulations, it can be arrived at (\cite{Kriv90})
$$
D_{n}(k) = \sumol_{1 \leq i_{1} \leq \cdots \leq i_{n} \leq n}
\left(\prodo_{j=1}^{i_{1}} \tau_{0j} \otimes
\prodo_{j=i_{1}}^{i_{2}} \tau_{1j} \otimes \cdots \otimes
\prodo_{j=i_{n}}^{k} \tau_{nj}\right), \;\;\; k=1,2, \ldots
$$

\subsection{Tandem Queues with Finite Buffers}
Suppose now that the buffers of servers in the tandem system described above
have finite capacity. The feature of queueing systems with limited buffers is
that their servers may be blocked according to one of the blocking rules
\cite{Gree91}. In this paper we restrict our consideration to
{\it manufacturing} blocking which is most commonly encountered in practice.
Under this type of blocking, if the $i$th server upon completion of a service
sees the buffer of the $(i+1)$st server full, it cannot be unoccupied and has
to be busy until the $(i+1)$st server completes its current service to provide
a free space in its buffer.

Consider a queueing system with $ n $ servers in tandem, and assume the
buffer at the $i$th server, $ i=2, \ldots, n$, to be of the capacity
$ b_{i}$, $ 0 \leq b_{i} < \infty $. We suppose that the buffer of the
$1$st server, as the input buffer of the system, is infinite. Since the
customers leave the system upon their service completion at the $n$th server,
this server cannot be blocked.

It is not difficult to understand that the $k$th completion time at the $i$th
server, $ i=1, \ldots, n-1 $, can be represented in usual form by the
recursive equation \cite{Gree91,Kriv90}
$$
D_{i}(k) = \max(\max(D_{i-1}(k), D_{i}(k-1)) + \tau_{ik}, D_{i+1}(k-b_{i+1}-1)).
$$
Using max--algebra notations, the complete set of linear equations describing
the finite buffers tandem queueing system with manufacturing blocking is
written as
\begin{eqnarray}
D_{0}(k) & = & \tau_{0k} \otimes D_{0}(k-1), \label{set3-1} \\
D_{i}(k) & = & \tau_{ik} \otimes (D_{i-1}(k) \oplus D_{i}(k-1))
                                             \oplus D_{i+1}(k-b_{i+1}-1),\nonumber
\\                                             
& & i=1, \ldots, n-1, \label{set3-2} \\
D_{n}(k) & = & \tau_{nk} \otimes (D_{n-1}(k) \oplus D_{n}(k-1)).\label{set3-3}
\end{eqnarray}

Although it is evident that we are dealing here with a linear model once
again, handling the model in its general form (\ref{set3-1}-\ref{set3-3})
requires rather cumbersome algebraic manipulations. Therefore, let us consider
more thoroughly a simple example of a system with $ n=2 $, $ b_{2}=0 $.
The equations (\ref{set3-1}-\ref{set3-3}) are reduced to
\begin{eqnarray}
D_{0}(k) & = & \tau_{0k} \otimes D_{0}(k-1), \label{set31-1} \\
D_{1}(k) & = & \tau_{1k} \otimes (D_{0}(k) \oplus D_{1}(k-1))
                                      \oplus D_{2}(k-1), \label{set31-2}\\
D_{2}(k) & = & \tau_{2k} \otimes (D_{1}(k) \oplus D_{2}(k-1)). \label{set31-3}
\end{eqnarray}
Going to matrix notations, we arrive at the linear equation
$$
\mbox{\boldmath $D$}(k) = \tilde{T}_{k} \otimes \mbox{\boldmath $D$}(k-1),
$$
with
$$
\tilde{T}_{k} = \left( \begin{array}{rrc}
         	        \tau_{0k} & \varepsilon & \varepsilon \\
		        \tau_{1k} \otimes \tau_{0k} & \tau_{1k} & e \\
		        \tau_{2k} \otimes \tau_{1k} \otimes \tau_{0k} &
                            \tau_{2k} \otimes \tau_{1k} & \tau_{2k}
	               \end{array} \right).
$$
The above representation of the transition matrix for the system with two
servers is easily extended to the case of the system with $ n $ servers
and $ b_{i}=0 $, $ i=2, \ldots, n$
$$
\tilde{T}_{k} = \left( \begin{array}{rrrcc}
	       \tau_{0k} & \varepsilon & \varepsilon & \ldots & \varepsilon \\
	       \tau_{1k} \otimes \tau_{0k} & \tau_{1k} &
                                                e &     & \varepsilon \\
	       \vdots & \vdots &     & \ddots &    \\
	       \tau_{n-1k} \otimes \cdots \otimes \tau_{0k} &
	       \tau_{n-1k} \otimes \cdots \otimes \tau_{1k} &
               \tau_{n-1k} \otimes \cdots \otimes \tau_{2k} &    & e \\
	       \tau_{nk} \otimes \cdots \otimes \tau_{0k} &
	       \tau_{nk} \otimes \cdots \otimes \tau_{1k} &
               \tau_{nk} \otimes \cdots \otimes \tau_{2k} &
                                                        \ldots & \tau_{nk}
	       \end{array} \right).
$$
Note that the matrices $ \tilde{T}_{k} $ and $ T_{k} $ differ only in
elements of the upper diagonal adjacent to the main diagonal. In
$ \tilde{T}_{k} $ these elements become equal to $ e $, excluding the
one of row $ 0 $ which remains equaled $ \varepsilon $.

Now we return to the example so as to present the solution of the recursive
equations (\ref{set31-1}-\ref{set31-3}). After usual algebraic manipulations
one can obtain
\begin{eqnarray*}
D_{1}(k) & = & \sumo_{i=1}^{k} \left( \prodo_{j=1}^{i} \tau_{0j} \otimes \tau_{1i}
\otimes \prodo_{j=i+1}^{k} (\tau_{1j} \oplus \tau_{2j-1}) \right),
\\
D_{2}(k) & = & \tau_{2k} \otimes D_{1}(k).
\end{eqnarray*}

\subsection{Closed Systems of $ G/G/1 $ Queues}
Consider a closed tandem system of $\,n\,$ queues with infinite buffers. We
assume that the customers after their service completion at the $n$th server
return to the $1$st server for a new cycle of service. There are the finite
number of customers circulating through the system, at the initial time all
the customers are placed in the buffer of the $1$st server.

Let us denote the number of customers in the system by $ c$. With the
condition $ D_{n}(k) = \varepsilon $ for all $ k < 0$, we may  modify
(\ref{set2-1}-\ref{set2-2}) to write the set of equations for the closed
system without specifying $ D_{0}(k)$, in the form
\begin{eqnarray*}
D_{1}(k) & = & \tau_{1k} \otimes (D_{n}(k-c) \oplus D_{1}(k-1)), \\
D_{i}(k) & = & \tau_{ik} \otimes (D_{i-1}(k) \oplus D_{i}(k-1)), \;\;\;
					i=2, \ldots, n.
\end{eqnarray*}

To produce a matrix representation, we now define the vector of departure
times as $ \mbox{\boldmath $D$}(k) = (D_{1}(k), \ldots, D_{n}(k))^{T} $.
The vector equation associated with the closed tandem queueing system with
$ c $ customers is represented as
$$
\mbox{\boldmath $D$}(k) = R_{k} \otimes \mbox{\boldmath $D$}(k-1) \oplus
                          S_{k} \otimes \mbox{\boldmath $D$}(k-c),
$$
where
\begin{eqnarray*}
R_{k} & = &\left( \begin{array}{rrcc}
		 \tau_{1k} & \varepsilon & \ldots & \varepsilon \\
		 \tau_{2k} \otimes \tau_{1k} & \tau_{2k} &    &
                 \varepsilon \\
		 \vdots & \vdots & \ddots &    \\
		 \tau_{nk} \otimes \cdots \otimes \tau_{1k} &
		 \tau_{nk} \otimes \cdots \otimes \tau_{2k} &
                                                         \ldots & \tau_{nk}
	       \end{array} \right),
\\	       
S_{k} & = &\left( \begin{array}{cccr}
		 \varepsilon & \ldots & \varepsilon & \tau_{1k} \\
		 \varepsilon & \ldots & \varepsilon &
                 \tau_{2k} \otimes \tau_{1k} \\
		 \vdots &        &  \vdots & \vdots \\
		 \varepsilon & \ldots & \varepsilon & \tau_{nk} \otimes \cdots
                                                        \otimes \tau_{1k}
	       \end{array} \right).
\end{eqnarray*}

In conclusion, consider an example of a closed system with $ n=2 $,
$ c=2 $. We have the following linear equations representing this system
\begin{eqnarray*}
D_{1}(k) & = & \tau_{1k} \otimes (D_{1}(k-1) \oplus D_{2}(k-2)), \\
D_{2}(k) & = & \tau_{2k} \otimes (D_{1}(k) \oplus D_{2}(k-1)).
\end{eqnarray*}
Traditional methods of solving linear recursions give the solution
\begin{eqnarray*}
D_{1}(k) & = & \tau_{11} \otimes
\prodo_{i=1}^{k-2} (\tau_{1i+1} \oplus \tau_{2i}) \otimes \tau_{1k},
\\
D_{2}(k) & = & \tau_{11} \otimes
\prodo_{i=1}^{k-1} (\tau_{1i+1} \oplus \tau_{2i}) \otimes \tau_{2k}.
\end{eqnarray*}

\bibliographystyle{utphys}

\bibliography{Using_max-algebra_linear_models_in_the_representation_of_queueing_systems}

\end{document}